\documentclass[12pt]{amsart}
\usepackage{amssymb,amscd}
\usepackage[all]{xy}

\theoremstyle{plain}
\newtheorem{theorem}{Theorem}[section]
\newtheorem{corollary}[theorem]{Corollary}
\newtheorem{lemma}[theorem]{Lemma}
\newtheorem{proposition}[theorem]{Proposition}

\theoremstyle{definition}

\newtheorem{remark}[theorem]{Remark}
\newtheorem*{question}{Question}

\theoremstyle{remark}

\newcommand{\abs}[1]{\lvert#1\rvert}
\newcommand{\norm}[1]{\lVert#1\rVert}
\newcommand{\bigabs}[1]{\bigl\lvert#1\bigr\rvert}
\newcommand{\bignorm}[1]{\bigl\lVert#1\bigr\rVert}

\newcommand{\Bignorm}[1]{\Bigl\lVert#1\Bigr\rVert}
\renewcommand{\le}{\leqslant}
\renewcommand{\ge}{\geqslant}

\newcommand{\term}[1]{{\textit{\textbf{#1}}}}

\newcommand{\hyph}{-\penalty0\hskip0pt\relax } 
\newcommand{\T}{Pe{\l}\-czy{\'n}\-ski Decomposition}

\def\Lpq{L(\ell_p,\ell_q)}
\def\Lpoq{L(\ell_p\oplus\ell_q)}
\def\iJ{\mathcal J}
\def\iK{\mathcal K}
\def\fss{\iJ^{\mathrm FSS}}
\def\Ipq{I_{p,q}}
\def\JIpq{\iJ^{\Ipq}}
\def\I2pq{I_{2,p,q}}
\def\JT{\iJ^T}
\def\J2{\iJ^{\ell_2}}
\def\Jr{\iJ^{\ell_r}}
\def\toplus{\textstyle\bigoplus}

\DeclareMathOperator{\supp}{supp}
\DeclareMathOperator{\Span}{span}

\DeclareMathOperator{\diag}{diag}
\DeclareMathOperator{\sign}{sign}
\DeclareMathOperator{\re}{Re}

\DeclareMathOperator{\card}{card}
\DeclareMathOperator{\codim}{codim}
\DeclareMathOperator{\rank}{rank}

\begin{document}
\baselineskip 18pt

\title{On norm closed ideals in $\Lpoq$.}
\author[Sari]{B.~Sari}
\address{B.~Sari,
          Department of Mathematics, University of North Texas,
          Denton, TX 76203-1430. USA.}
\email{bunyamin@unt.edu}
\author[Schlumprecht]{Th.~Schlumprecht}
\address{Th.~Schlumprecht, 
         Department of Mathematics, Texas A\&M University,
         College Station, TX 77843-3368. USA.}
\email{schlump@math.tamu.edu}
\author[Tomczak-Jaegermann]{N.~Tomczak-Jaegermann}
\address{N.~Tomczak-Jaegermann \and V.G.~Troitsky, 
         Department of Mathematical and Statistical Sciences,
         University of Alberta, Edmonton, AB, T6G\,2G1. Canada.}
\email{nicole@ellpspace.math.ualberta.ca}
\author[Troitsky]{V.G.~Troitsky} 
\email{vtroitsky@math.ualberta.ca}

\thanks{The first author was supported by the University of Alberta
  postdoctoral fellowship. The second author was supported by NSF. The
  third author holds the Canada Research Chair in Geometric Analysis.
  The fourth author was supported by the University of Alberta
  start-up grant. Most of the work on the paper was done during second
  author's visit to the University of Alberta in 2003 and during the
  {\it Workshop on linear analysis and probability} at Texas A$\&$M
  University in 2004.}  

\date{\today.}

\begin{abstract}
  It is well known that the only proper non-trivial norm-closed ideal
  in the algebra $L(X)$ for $X=\ell_p$ $(1\le p<\infty)$ or $X=c_0$ is
  the ideal of compact operators. The next natural question is to
  describe all closed ideals of $\Lpoq$ for $1\le p,q<\infty$, $p\neq
  q$, or, equivalently, the closed ideals in $\Lpq$ for $p<q$. This
  paper shows that for $1<p<2<q<\infty$ there are at least four
  distinct proper closed ideals in $\Lpq$, including one that has not
  been studied before. The proofs use various methods from Banach space
  theory.
\end{abstract}

\maketitle

\section{Introduction}

This paper is concerned with the structure of norm closed ideals of
the algebra $L(X)$ of all bounded linear operators on an
infinite-dimensional Banach space~$X$. The classical result
of~\cite{Calkin:41} asserts that the only proper non-trivial ideal of
$L(\ell_2)$ is the ideal of compact operators. The same was shown to
be true for $\ell_p$ $(1\le p<\infty)$ and $c_0$ in~\cite{Gohberg:60}.
It remains open if there are other Banach spaces with only one proper
non-trivial closed ideal.  The complete structure of closed ideals in
$L(X)$ was recently described in~\cite{Laustsen:04} for
$X=\bigl(\bigoplus_{n=1}^\infty\ell_2^n\bigr)_{c_0}$ and
in~\cite{Laustsen:p} for
$X=\bigl(\bigoplus_{n=1}^\infty\ell_2^n\bigr)_{\ell_1}$. In the both
cases, there are exactly two nested proper non-zero closed ideals.
Apart from those mentioned above, there are no other Banach spaces $X$
for which the structure of the closed ideals in $L(X)$ is completely
known.

This motivates the study of the next natural special case
$X=\ell_p\oplus\ell_q$ ($1\le p,q<\infty$, $p\neq q$), which is our
main interest here.  There were several results in this direction
proved in the 1970's concerning various special ideals or special
cases of $p$ and~$q$. We refer the reader to the book by Pietsch
\cite[Chapter~5]{Pietsch:78} for details. In particular,
\cite[Theorem~5.3.2]{Pietsch:78} asserts that $\Lpoq$ (with, say,
$p<q$) has exactly two proper maximal ideals (namely, the ideal of
operators which factor through $\ell_p$ and the ideals of operators
which factor through $\ell_q$), and establishes a one-to-one
correspondence between the non-maximal ideals in the algebra $\Lpoq$
and the closed ``ideals'' in $\Lpq$. Here an \term{ideal} in
$L(\ell_p,\ell_q)$ means a linear subspace $\iJ$ of $L(\ell_p,\ell_q)$
such that $ATB\in\iJ$ whenever $A\in L(\ell_q)$, $T\in\iJ$, and $B\in
L(\ell_p)$, and ``closed'' is always understood with respect to the
operator norm topology.  Consequently, the subject of the present
paper is a study of the structure of closed ideals in $\Lpq$ with $1\le
p<q<\infty$.

In this paper, we identify four distinct proper closed ideals in
$\Lpq$ when $1<p<2<q<\infty$ (however, some of the results remain
valid under weaker restrictions on $p$ and~$q$). Namely, the ideal of
all compact operators $\mathcal K$, the closed ideal $\JIpq$ generated
by the formal identity operator $\Ipq\colon\ell_p\to\ell_q$, the ideal
of all finitely strictly singular operators $\fss$, and the closure of
the ideal of all $\ell_2$-factorable operators $\J2$ (see
Section~\ref{notation} for appropriate definitions). Although
these ideals were recognized earlier, they were not known to be
distinct and proper except for special cases of $p$ and~$q$.  The
following diagram illustrates the relationship between these ideals.

\noindent
{\footnotesize
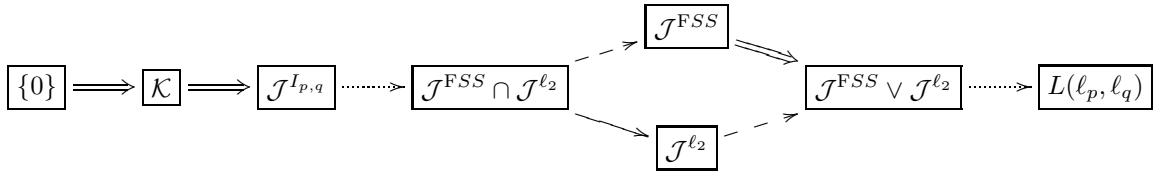
\xymatrix@R=0pt@C=23pt{    
& & & & \boxed{\fss} \ar@{=>}[dr] & &\\
\boxed{\{0\}} \ar@{=>}[r] & \boxed{\mathcal K} \ar@{=>}[r] & 
  \boxed{\JIpq} \ar@{.>}[r] &
  \boxed{\fss\cap\J2} \ar@{-->}[ur] \ar[dr] & &
  \boxed{\fss\vee\J2} \ar@{.>}[r] & \boxed{\Lpq}\\
& & & & \boxed{\J2} \ar@{-->}[ur] & &
}}\\
Here arrows stand for inclusions. A solid arrow ($\Rightarrow$ or
$\to$) between two ideals means that there are no other ideals
sitting properly between the two, while a double arrow comming out of
an ideal indicates the only immediate successor.  A
hyphenated arrow ($--\!\!\!>$) indicates a proper inclusion, while a dotted
one indicates that we do not know whether or not the inclusion is
proper. In particular, the closed ideals in $\Lpq$ are not totally
ordered.

The paper is organized as follows. In Section~\ref{sec:Ipq} we study
the ideal $\JIpq$ for $1\le p<q<\infty$.  In~\cite{Milman:70}, Milman
proved that $\JIpq$ is FSS, and, therefore, $\JIpq\subseteq\fss$.
Since $\JIpq$ is not compact, $\iK$ is properly contained in $\JIpq$.
We will show that every closed ideal that contains a non-compact
operator necessarily contains $\JIpq$, so that $\JIpq$ is the least
non-compact ideal.  In Section~\ref{sec:hilbert} we consider the ideal
$\J2$ when $1<p\le 2\le q<\infty$. We find a specific non-FSS operator
$T$ in $\J2$ such that the closed ideal $\JT$ generated by $T$
coincides with $\J2$. This implies, in particular, that $\fss$ is a
proper ideal. It should be noted here that Milman proved in
\cite{Milman:70} that $\fss$ is a proper ideal for special values of
$p$ and~$q$.  We also show that $\J2\subseteq\Jr$ for all $r$ between
$p$ and~$q$. We prove that every closed ideal of $\Lpq$ which contains
a non-FSS operator must also contain $\J2$. In Section~\ref{sec:walsh}
we consider the ``block Hadamard'' operator $U$ from $\ell_p$ to
$\ell_q$ for $p<2<q$.  We show that $U\notin\J2$, hence $\J2$ is a
proper ideal. Since, obviously, $\Ipq\in\J2$, it follows that
$\JIpq\subsetneq\iJ^U$.  We show in Section~\ref{Ufss} that $U$ is
FSS, hence $\JIpq\subsetneq\fss$.

We thank Gilles Pisier for suggesting to us the proof of
Theorem~\ref{dimension}.

\section{Notation and preliminaries}\label{notation}

Given two Banach spaces $X$ and $Y$, we write $L(X,Y)$ for the space
of all continuous linear operators from $X$ to $Y$, we write $L(X)$
for $L(X,X)$. A linear subspace $\iJ$ of $L(X,Y)$ is said to be an
\term{ideal} if $ATB\in\iJ$ whenever $A\in L(Y)$, $T\in\iJ$, and $B\in
L(X)$. By a \term{closed ideal} we mean an ideal closed in the
operator norm topology. We denote by $\iK$ the closed ideal of all
compact operators.

Throughout this paper, $p$ and $q$ always satisfy $1\le p<q<\infty$.
We denote by $p'$ the conjugate of~$p$, that is,
$\frac{1}{p}+\frac{1}{p'}=1$.  It is well known (see,
e.g.,~\cite{Caradus:74}) that $\iK$ is contained in every closed ideal
of $\Lpq$. If $Z$ is a Banach space, we denote by $\iJ^Z$ the closure
of the set of all the operators in $\Lpq$ that factor through $Z$. It
can be easily verified that if $Z$ is isomorphic to $Z\oplus Z$ then
$\iJ^Z$ is a subspace, hence an ideal. For $S\in\Lpq$ we denote by
$\iJ^S$ the closed ideal in $\Lpq$ generated by $S$, that is, the
smallest closed ideal containing~$S$.  It is easy to see that $\iJ^S$
consists of the operators that can be approximated in norm by
operators of the form $\sum_{i=1}^nA_iSB_i$, where $A_i\in L(\ell_q)$
and $B_i\in L(\ell_p)$ for $i=1,\dots,n$.  If $A$ is an $n\times n$
scalar matrix, we write $\norm{A}_{p,q}$ for the norm of $A$ as an
operator from $\ell_p^n$ to $\ell_q^n$.
  
It is known that every operator in $\Lpq$ is strictly singular, see,
e.g.,~\cite{Lindenstrauss:77}. We call an operator $S\colon X\to Y$
\term{finitely strictly singular} or \term{FSS} if for every
$\varepsilon>0$ there exists $n\in\mathbb N$ such that
$\inf\limits_{x\in E,\;\norm{x}=1}\norm{Sx}<\varepsilon$ for every
$n$\hyph dimensional subspace $E$ of~$X$. It can be easily verified
(see~\cite{Mascioni:94}) that $S$ is FSS if and only if every
ultrapower of $S$ is strictly singular. It follows immediately that
the set of all FSS operators from $X$ to $Y$ is a closed ideal. Denote
by $\fss$ the ideal of all FSS operators in $\Lpq$.

We denote by $(e_i)$ and $(f_i)$ the standard bases of $\ell_p$ and
$\ell_q$ respectively, and we denote their coordinate functionals by
$(e_i^*)$ and $(f_i^*)$. If $(x_n)$ is a sequence in a Banach space,
we write $[x_n]$ for its closed linear span. A sequence $(x_n)$ in a
Banach space is \term{semi-normalized} if $\inf_n\norm{x_n}>0$ and
$\sup_n\norm{x_n}<\infty$.

The following standard lemma follows immediately from
Propositions~1.a.12 and~2.a.1 of~\cite{Lindenstrauss:77}.

\begin{lemma}\label{standard}
  If $X=\ell_p$ ($1\le p<\infty$) or $c_0$ and $(x_n)$ is a
  semi-normalized sequence in $X$ which converges to zero
  coordinate-wise (that is, $e_i^*(x_n)\to 0$ in $n$ for every $i$),
  then there is a subsequence $(x_{n_i})$ equivalent to $(e_i)$, and
  $[x_{n_i}]$ is complemented in $X$.
\end{lemma}

\begin{remark}\label{block-diag}
  Suppose that $1\le p\le q<\infty$ and $T\in\Lpq$. We say that $T$ is
  \term{block-diagonal} if $T=\bigoplus_{n=1}^\infty T_n$, where
  $T_n\colon\ell_p^{m_n}\to\ell_q^{m_n}$. Equivalently, there exists a
  strictly increasing sequence of integers $(k_n)$ such that
  $T=\sum_{n=1}^\infty P_nTQ_n$, where $Q_n$ and $P_n$ are the
  canonical projections from $\ell_p$ and $\ell_q$ to the
  finite-dimensional subspaces spanned by
  $e_{k_n+1},\dots,e_{k_{n+1}}$ and $f_{k_n+1},\dots,f_{k_{n+1}}$
  respectively. Note that $m_n=k_{n+1}-k_n$ and $T_n$ can be
  identified with $P_nTQ_n$. It can be easily verified that if $p\le
  q$ then $\norm{T}=\sup_n\norm{T_n}$.  Indeed,
  $\norm{T_n}=\norm{P_nTQ_n}\le\norm{T}$ as $P_n$ and $Q_n$ are
  contractions. On the other hand,
  \begin{multline*}
   \norm{Tx}=
   \Bigl(\sum_{n=1}^\infty\norm{P_nTQ_nx}^q\Bigr)^{\frac{1}{q}}\le
   \Bigl(\sup\limits_{n}\norm{P_nTQ_n}\Bigr)
   \Bigl(\sum_{n=1}^\infty\norm{Q_nx}^q\Bigr)^{\frac{1}{q}}\\\le
   \bigl(\sup\limits_{n}\norm{T_n}\bigr)
   \Bigl(\sum_{n=1}^\infty\norm{Q_nx}^p\Bigr)^{\frac{1}{p}}=
   \bigl(\sup\limits_{n}\norm{T_n}\bigr)\norm{x}.
  \end{multline*}
\end{remark}

\begin{remark}\label{block-convex}
  Suppose that $R\in\Lpq$ for $1\le p\le q<\infty$, and $T$ is a
  block-diagonal submatrix of $R$, that is, $T=\sum_{n=1}^\infty
  P_nRQ_n$, where $(P_n)$ and $(Q_n)$ are as in
  Remark~\ref{block-diag}. Then $T$ can be written as a convex
  combination of operators of the form $URV$, where $U$ and $V$ are
  isometries. See Proposition~1.c.8 of~\cite{Lindenstrauss:77} and
  Remark~1 following it for the construction.
\end{remark}

\section{The formal identity operator $\Ipq$}
\label{sec:Ipq}

In this section we consider the formal identity operator
$\Ipq\colon\ell_p\to\ell_q$ for $1\le p<q<\infty$.  Clearly, $\Ipq$ is
not compact, so that $\mathcal K\subsetneq\JIpq$.  First, we show that
$\JIpq$ is contained in every closed ideal of $\Lpq$ except $\iK$.
This result is probably known to specialists, but we provide a short
proof for completeness.

\begin{proposition}\label{Ipq}
  Let $1\le p<q<\infty$. If $\iJ$ is any ideal in $\Lpq$ containing a
  non-compact operator, then $\Ipq\in\iJ$.
\end{proposition}

\begin{proof}
  Assume that $\iJ$ contains a non-compact operator~$T$.  There exists
  a normalized sequence $(x_n)$ in $\ell_p$ such that $(Tx_n)$ has no
  convergent subsequences. By passing subsequences and using a
  standard diagonalization argument, we can assume that $(x_n)$ and
  $(Tx_n)$ converge coordinate-wise. Let $y_n=x_n-x_{n-1}$, then
  $(y_n)$ and $(Ty_n)$ converge coordinate-wise to zero. Since
  $(Tx_n)$ has no convergent subsequences, we can assume (by passing
  to a further subsequence if necessary) that $(Ty_n)$ is
  semi-normalized. It follows that $(y_n)$ is also semi-normalized.
  Using Lemma~\ref{standard} twice, we can assume (by passing to a
  subsequence) that $(y_n)$ is equivalent to $(e_i)$, $(Ty_n)$ is
  equivalent to $(f_i)$, and $[Ty_n]$ is complemented in $\ell_q$.

  Let $B\colon \ell_p\to[y_n]$ be an isomorphism given by $Be_n=y_n$,
  and let $A\colon[Ty_n]\to \ell_q$ be an isomorphism given by
  $A(Ty_n)=f_n$. Since $[Ty_n]$ is complemented, $A$ can be extended
  to an operator on all of $\ell_q$. Thus, we can view $B$ and $A$ as
  elements of $L(\ell_p)$ and $L(\ell_q)$ respectively. Observe that
  $ATBe_n=f_n$ for each $n$, hence $ATB=\Ipq$. It follows that
  $\Ipq\in\iJ$.
\end{proof}

\begin{corollary}
  If a closed ideal of $\Lpq$ contains a non-compact operator, then it
  contains $\JIpq$. 
\end{corollary}

The following result was proved in~\cite{Milman:70}.  For the
convenience of the reader we provide a short proof of it.

\begin{proposition}\label{Ipq-fss}
  Suppose that $1\le p<q<\infty$. The formal identity operator $\Ipq$
  is FSS.
\end{proposition}

We will deduce this proposition from the following lemma, which
appeared in~\cite{Milman:70}.

\begin{lemma}\label{flat-max}
  If $E$ is an $n$-dimensional subspace of $c_0$ then there exists $x\in
  E$ such that $x$ attains its sup-norm at at least $n$ coordinates.
\end{lemma}

\begin{proof}
  The proof is by induction. The statement is trivial for $n=1$.
  Suppose that it is true for $n$, take any subspace $E$ of $c_0$ of
  dimension $n+1$. By induction hypothesis, there exists $x\in E$ such
  that
  \begin{equation}\label{eq:max}
    \delta:=\norm{x}_\infty=\abs{x_{i_1}}=\dots=\abs{x_{i_n}}
  \end{equation}
  for a set of distinct indices $I=\{i_1,\dots,i_n\}$. Suppose that
  $\abs{x_i}<\delta$ for all $i\notin I$ (otherwise we are done).  Let
  $Y$ be the subspace of $c_0$ consisting of all the sequences that
  vanish at $i_1,\dots,i_n$. Since $Y$ has co-dimension $n$, it
  follows that $Y\cap E\neq\{0\}$. Pick a non-zero $y\in Y\cap E$.  We
  claim that for some $s>0$ the sequence $x+s y$ attains its sup-norm
  at at least $n+1$ coordinates.  Indeed, $\abs{x_i+ty_i}=\delta$ for
  all $i\in I$ and $t\ge 0$. Consider the function
  $$f(t)=\max\limits_{j\notin I}\abs{x_j+t y_j}.$$
  Clearly, $f$ is continuous, $f(0)<\delta$, and
  $\lim_{t\to+\infty}f(t)=+\infty$. It follows that $f(s)=\delta$ for
  some $s>0$. Then $\abs{x_i+s y_i}=\norm{x+sy}_\infty=\delta$ for
  some $i\notin I$.
\end{proof}

\begin{proof}[Proof of Proposition~\ref{Ipq-fss}]
  Given $\varepsilon>0$, let $n\in\mathbb N$ such that
  $n^{\frac{1}{q}-\frac{1}{p}}<\varepsilon$. Suppose that $E$ is a
  subspace of $\ell_p$ with $\dim E=n$.  By Lemma~\ref{flat-max} there
  exists $x\in E$ and indices $i_1,\dots,i_n$
  satisfying~(\ref{eq:max}). Without loss of generality,
  $\norm{x}_p=1$. It follows that $1=\norm{x}_p^p\ge n\delta^p$, so
  that $\delta\le n^{-\frac{1}{p}}$. Then
  $$\norm{x}_q^q\le\norm{x}_\infty^{q-p}\norm{x}_p^p=\delta^{q-p}\le
  n^{-\frac{1}{p}(q-p)},$$
  so that $\norm{x}_q\le
  n^{\frac{1}{q}-\frac{1}{p}}<\varepsilon$.  It follows that $\Ipq$ is
  FSS.
\end{proof}

\begin{corollary}
  Let $1\le p<q<\infty$. The ideal $\mathcal K$ is a proper subset of
  $\fss$.
\end{corollary}


\section{Operators factorable through $\ell_2$}
\label{sec:hilbert}

In this section we consider the ideal $\J2$ for $1<p<2<q$. Using
Pe{\l}czy{\'n}ski's decomposition, we will construct an operator
$T\colon\ell_p\to\ell_q$ such that $\J2=\iJ^T$. That is, the closure
of the ideal of all $\ell_2$-factorable operators is exactly the
closed ideal generated by~$T$.  Furthermore, we show that $T$ fails to
be FSS, hence the ideal $\fss$ is proper. It will be obvious from the
definition of $T$ that $T$ factors through $\ell_r$ whenever $p\le
r\le q$, so it follows that $\J2\subseteq\Jr$. We also show that $T$
factors through every non-FSS operator. It follows that any closed
ideal containing a non-FSS operator necessarily contains $\J2$.

To construct~$T$, recall that it follows from Pe{\l}czy{\'n}ski's
Decomposition Theorem that for every $1<r<\infty$, $\ell_r$ is isomorphic
to $\bigl(\bigoplus_{n=1}^\infty \ell_2^n\bigr)_r$, the
$\ell_r$-direct sum of $\ell_2^n$'s (see
\cite[p.~73]{Lindenstrauss:77}).  Let $p<q$, put
$U\colon\ell_p\to\bigl(\bigoplus_{n=1}^\infty \ell_2^n\bigr)_p$ and
$V\colon\bigl(\bigoplus_{n=1}^\infty \ell_2^n\bigr)_q\to\ell_q$ be two
such isomorphisms. By
$\I2pq\colon\bigl(\bigoplus_{n=1}^\infty\ell_2^n\bigr)_p
\to\bigl(\bigoplus_{n=1}^\infty\ell_2^n\bigr)_q$ we denote the formal
identity operator, that is, just the change of the norm on the direct
sum. Then let $T=V\I2pq U$, that is,
\begin{equation}\label{T}
  T\colon\ell_p\xrightarrow{U}\Bigl(\bigoplus\limits_{n=1}^\infty \ell_2^n\Bigr)_p
   \xrightarrow{\I2pq}\Bigl(\bigoplus\limits_{n=1}^\infty \ell_2^n\Bigr)_q 
  \xrightarrow{V}\ell_q.
\end{equation}
We will call $T$ a \term{\T} operator.

\begin{remark}\label{T-like}
  Note that $T$ is not unique, it is defined up to the isomorphisms
  $U$ and $V$, so that we have actually constructed a class of
  operators.  It is clear, however, that any two \T{} operators factor
  through each other.  Moreover, one can easily verify that if in the
  preceding construction we ``skip'' some of the blocks, that is, if
  we consider $\bigl(\bigoplus_{n=1}^\infty \ell_2^{k_n}\bigr)$ for
  some increasing sequence of indices $k_n$ then the resulting
  operator $T'$ obviously factors though~$T$. Conversely, $T$ factors
  through~$T'$ because $\ell_2^n$ is a complemented subspace of
  $\ell_2^{k_n}$.
\end{remark}

Furthermore, let $E_n=U^{-1}(\ell_2^n)\subset\ell_p$ be the pre-image
of the $n$-th block of $(\bigoplus \ell_2^n)_p$. Similarly, put
$F_n=V(\ell_2^n)\subset\ell_q$.  Then
$d(E_n,\ell_2^n)\le\norm{U}\cdot\norm{U^{-1}}$ and
$d(F_n,\ell_2^n)\le\norm{V}\cdot\norm{V^{-1}}$, where $d(X,Y)$ stands
for the Banach-Mazur distance between $X$ and $Y$.  Hence, $(E_n)$ and
$(F_n)$ are sequences of uniformly Euclidean subspaces of $\ell_p$ and
$\ell_q$ respectively. Note that $T(E_n)=F_n$, so that $T$ fixes
copies of $\ell_2^n$ for all $n\in\mathbb N$. This immediately implies
the following result.

\begin{proposition}
  For $p<q$, every \T{} operator fails to be FSS.
\end{proposition}

\begin{corollary}
  For $p<q$, the ideal $\fss$ is proper.
\end{corollary}

Our next goal is to show that if $1<p\le 2\le q<\infty$ then
$\iJ^T=\J2$.  We will make use of the concept of $\ell_2$-factorable
norm $\gamma_2$.  Recall that if $S\in L(X,Y)$ then
$\gamma_2(S)=\inf\norm{A}\norm{B}$, where the infimum is taken over
all factorizations $S=AB$ where $B\colon X\to\ell_2$ and
$A\colon\ell_2\to Y$. It is known that $\gamma_2$ is a norm on the
ideal of all $\ell_2$-factorable operators, and
$\gamma_2(ASB)\le\norm{A}\gamma_2(S)\norm{B}$ whenever
$X\xrightarrow{B}X\xrightarrow{S}Y\xrightarrow{A}Y$.
See~\cite{Tomczak:89,Diestel:95} for more information on $\gamma_2$.

\begin{lemma}\label{skipping-blocks}
  Suppose that $R\in\Lpq$, $1<p\le q<\infty$, and $\varepsilon>0$.
  \begin{enumerate}
  \item\label{op-norm}
    There exist two block-diagonal operators $V,W\in\Lpq$ such that
    $\norm{W}\le\norm{R}$, $\norm{V}\le2\norm{R}+\varepsilon$, and
    $\bignorm{R-(W+V)}<\varepsilon$.
  \item\label{gamma2} Suppose that, in addition, $R$ is
    $\ell_2$-factorable. Then $V$ and $W$ can be chosen to be
    $\ell_2$-factorable, and $\gamma_2(W)\le\gamma_2(R)+\varepsilon$,
    $\gamma_2(V)\le2\gamma_2(R)+\varepsilon$, and
    $\gamma_2\bigl(R-(W+V)\bigr)<\varepsilon$.
  \end{enumerate}
\end{lemma}

\begin{proof}
  Let $r_{i,j}$ stand for the $(i,j)$-th entry of the matrix of $R$,
  that is, $r_{i,j}=f_i^*(Re_j)$. One can approximate $R$ with a
  matrix $S$ with finitely many entries in every row and every column.
  That is, there exists an operator $S=(s_{i,j})$ and two strictly
  increasing sequences $(M_j)$, $(N_i)$ of positive integers
  such that $\norm{R-S}<\varepsilon$ and
  \begin{displaymath}
    s_{i,j}=
    \begin{cases}
      r_{i,j}&\mbox{if } i\le M_j\mbox{ and }j\le N_i;\\
            0&\mbox{otherwise.}
    \end{cases}
  \end{displaymath}
  Let $\Gamma$ be the subset of $\mathbb N\times\mathbb N$ consisting of
  all the pairs of indices corresponding to the ``non-trivial'' part
  of~$S$, namely,
  \begin{equation*}
    (i,j)\in\Gamma \mbox{ iff } i\le M_j \mbox{ and } j\le N_i.
  \end{equation*}

  We will define two strictly increasing sequences $(k_n)$ and $(l_n)$
  of positive integers, such that $\Gamma$ is contained in the union of
  two block-diagonal sets
  $\Delta=\bigcup_{n=1}^\infty\Delta_n$ and 
  $\Lambda=\bigcup_{n=1}^\infty\Lambda_n$ where
  \begin{eqnarray*}
    \Delta_n &=& \bigl\{(i,j)\in\Gamma\mid k_{n-1}<i,j\le k_n\bigr\}
    \mbox{ and }\\
    \Lambda_n &=& \bigl\{(i,j)\in\Gamma\mid l_{n-1}<i,j\le l_n\bigr\}.
  \end{eqnarray*}
  We define the sequences $(k_n)$ and $(l_n)$ by an interlaced induction. Put
  $k_0=0$, $l_0=1$. For $n\ge 0$ we let
  \begin{equation*}
  k_{n+1}=\max\{M_{l_n},N_{l_n}\}\quad\mbox{ and }\quad
  l_{n+1}=\max\{M_{k_{n+1}},N_{k_{n+1}}\}.
  \end{equation*}
  Clearly, $(k_n)$ and $(l_n)$ are strictly increasing. Next, we show
  that $\Gamma\subseteq\Delta\cup\Lambda$. Let $(i,j)\in\Gamma$. There
  exists $n$ such that $l_n<\max\{i,j\}\le l_{n+1}$.  If
  $l_n<\min\{i,j\}$, then $l_n<i,j\le l_{n+1}$, so that
  $(i,j)\in\Lambda$.  Suppose now that $\min\{i,j\}\le l_n$. Then either
  $i$ or $j$ is less than or equal to~$l_n$, while the other is
  greater than~$l_n$.  Say, $i\le l_n$ and $j>l_n$. It follows that
  \begin{displaymath}
    i\le l_n\le N_{l_n}\le k_{n+1}\quad\mbox{ and }\quad
    j>l_n\ge N_{k_n}\ge k_n.
  \end{displaymath}
  Therefore $j\le N_i\le N_{l_n}\le k_{n+1}$.
  Also, $N_i\ge j>l_n\ge N_{k_n}$ yields $i>k_n$.
  Hence, $k_n<i,j\le k_{n+1}]$, so that $(i,j)\in\Delta$.

  Let $W=(w_{i,j})$ be the operator defined by
  \begin{equation}\label{W}
    w_{i,j}=
    \begin{cases}
      s_{i,j} & \mbox{if }(i,j)\in\Delta,\mbox{ and}\\
      0 & \mbox{otherwise}.
    \end{cases}
  \end{equation}
  Put $V=S-W$. Then the non-zero entries of $W$ and $V$ are located in
  $\Delta$ and $\Lambda$ respectively, so that $W$ and $V$ are
  block-diagonal. By the definition of $S$ we have
  $\bignorm{R-(W+V)}<\varepsilon$.


  Since $W$ is a block-diagonal part of $R$ Remark~\ref{block-convex}
  yields that $\norm{W}\le\norm{T}$.  Finally, it follows from $V=S-W$
  that $\norm{V}\le2\norm{R}+\varepsilon$.
  
  If $R$ is $\ell_2$-factorable, then we can choose $S$ with finitely
  many entries in each row and column such that $S$ is also
  $\ell_2$-factorable and $\gamma_2(R-S)<\varepsilon$. Indeed, let
  $R=R_1R_2$ be a factorization of $R$ through $\ell_2$. Approximate
  $R_1$ and $R_2$ in norm by $S_1$ and $S_2$ respectively, such that
  $S_1$ and $S_2$ have finitely many entries in every row and column.
  Put $S=S_1S_2$, then $S$ is as claimed. We use triangle inequality
  to show that $\gamma_2(R-S)<\varepsilon$ when $\norm{R_1-S_1}$ and
  $\norm{R_2-S_2}$ are sufficiently small.

  
  Define $W$ as a block-diagonal part of $S$ using~(\ref{W}).  It
  follows from Remark~\ref{block-convex} that
  $\gamma_2(W)\le\gamma_2(S)\le\gamma_2(R)+\varepsilon$.  Put $V=S-W$,
  then $\gamma_2(V)=\gamma_2(S-W)\le2\gamma_2(R)+\varepsilon$. In
  particular, $W$ and $V$ are $\ell_2$-factorable.
\end{proof}

\begin{remark}
  In a similar fashion one can show that every operator between two
  Banach spaces with shrinking unconditional bases can be approximated
  by a sum of two block-diagonal operators.
\end{remark}

\begin{theorem}\label{JT-is-J2}
  If $1<p\le 2\le q$ and $T$ is a \T{} operator, then $\JT=\J2$.
\end{theorem}

\begin{proof}
  Observe that $\I2pq$, being the formal identity from
  $\bigl(\bigoplus_{n=1}^\infty\ell_2^n\bigr)_p$ to
  $\bigl(\bigoplus_{n=1}^\infty\ell_2^n\bigr)_q$, factors through
  $\bigl(\bigoplus_{n=1}^\infty\ell_2^n\bigr)_2=\ell_2$. It follows
  that $T$ factors through $\ell_2$ and, therefore, $\JT\subseteq\J2$.
 
  We show that $\J2\subseteq\JT$. Clearly, it suffices to show that
  every $\ell_2$-factorable operator belongs to~$\JT$.  In view of
  Lemma~\ref{skipping-blocks}(\ref{gamma2}), it suffices to show this
  for block-diagonal operators.  Let $W$ be an $\ell_2$\hyph
  factorable block-diagonal operator. Then we can write
  $W=\bigoplus_{n=1}^\infty A_nB_n$, where
  $B_n\colon\ell_p^{k_n}\to\ell_2^{k_n}$ and
  $A_n\colon\ell_2^{k_n}\to\ell_q^{k_n}$ such that $\sup_n\norm{A_n}$
  and $\sup_n\norm{B_n}$ are finite. By merging consequtive blocks if
  necessary, we can assume without loss of generality that $(k_n)$ is
  strictly increasing. Observe that the operators
  \begin{eqnarray*}
    B&=&\bigoplus\limits_{n=1}^\infty B_n\colon
      \Bigl(\bigoplus\limits_{n=1}^\infty\ell_p^{k_n}\Bigr)_p\to
      \Bigl(\bigoplus\limits_{n=1}^\infty\ell_2^{k_n}\Bigr)_p
      \quad\mbox{and}\\
    A&=&\bigoplus\limits_{n=1}^\infty A_n\colon
      \Bigl(\bigoplus\limits_{n=1}^\infty\ell_2^{k_n}\Bigr)_q\to
      \Bigl(\bigoplus\limits_{n=1}^\infty\ell_q^{k_n}\Bigr)_q
  \end{eqnarray*}
  are bounded, and $W=AI_0 B$, where $I_0$ is the formal identity
  between $\bigl(\bigoplus_{n=1}^\infty\ell_2^{k_n}\bigr)_p$ and
  $\bigl(\bigoplus_{n=1}^\infty\ell_q^{k_n}\bigr)_q$. Thus, $W$ factors
  through~$I_0$. It follows from Remark~\ref{T-like} that $I_0$
  factors through~$T$. Hence, $W$ factors through~$T$.
\end{proof}

\begin{remark}
  Actually, we proved that every operator in $\J2$ can be approximated
  by sums of two T\hyph factorable operators.
\end{remark}

\begin{remark}\label{Iprq}
  Suppose that $p<r<q$. Then $\I2pq$ in~(\ref{T}) factors through
  $\bigl(\bigoplus_{n=1}^\infty\ell_2^n\bigr)_r$, which is isomorphic
  to~$\ell_r$. It follows that $T$ factors through $\ell_r$. Then
  Theorem~\ref{JT-is-J2} implies that $\J2\subseteq\Jr$.
\end{remark}

Next, we show that if $p<2<q$ then $\J2$ is the least closed ideal
beyond $\fss$, that is, every closed ideal that contains a non-FSS
operator also contains $\J2$. For the proof we need the following
well-known fact, which can be viewed, for example, as a special case
of results in~\cite{Figiel:79}.

\begin{theorem}\label{K-convex}
  For every $1<r<\infty$ there exists $K>0$ such that for all
  $n\in\mathbb N$ there exists $N\in\mathbb N$
  such that every $N$\hyph dimensional subspace $F\subset\ell_r$
  contains an $n$\hyph dimensional subspace $E$ which is $K$\hyph
  complemented in $\ell_r$ and 2\hyph isomorphic to~$\ell_2^n$.
\end{theorem}

We will also routinely use the following observation.

\begin{remark}\label{sub-unif}
  Suppose that $(E_n)$ is a sequence of subspaces of a Banach space
  $X$ such that for each $n$ we have $\dim E_n=n$ and $E_n$'s are
  uniformly Euclidean and uniformly complemented in~$X$. That is,
  there exist sequences $(P_n)$ and $(V_n)$ and a constant $C>0$ such
  that $P_n$ is a projection from $X$ onto $E_n$ with $\norm{P_n}<C$,
  and $V_n\colon E_n\to\ell_2^n$ is an isomorphism with
  $\norm{V_n}\cdot\norm{V_n^{-1}}\le C$ for every~$n$. For a
  subsequence $(E_{k_n})$, let $G_n$ be a subspace of $E_{k_n}$ with
  $\dim G_n=n$ for every~$n$. It is easy to see that $G_n$'s are still
  uniformly Euclidean and uniformly complemented in~$X$.
\end{remark}

For $x\in\ell_r$ we write $\supp x=\{i\in\mathbb N\mid x_i\ne 0\}$.
For $A\subseteq\ell_r$ put $\supp A=\cup_{x\in A}\supp x$.

\begin{theorem}
  Let $1<p\le 2\le q<\infty$. If $R\in\Lpq$ is not FSS, then every
  \T{} operator factors through~$R$.
\end{theorem}

\begin{proof}
  Since $R$ is not FSS, there exist a constant $C>0$ and a sequence
  $(E_n)$ of subspaces of $\ell_p$ such that $\dim E_n=n$, and
  $R_{|E_n}$ is invertible with $\norm{(R_{|E_n})^{-1}}\le C$. We can
  assume, in addition, that $\supp E_n$ is finite by truncating all
  the vectors in a basis of $E_n$ sufficiently far (and adjusting $C$
  if necessary).  Let $F_n=R(E_n)$. By Theorem~\ref{K-convex} and
  Remark~\ref{sub-unif}, we can also assume that the sequences $(E_n)$
  and $(F_n)$ are $C$-complemented in $\ell_p$ and $\ell_q$
  respectively, and $C$-isomorphic to~$\ell_2^n$. Let
  $Q_n\colon\ell_q\to F_n$ be a projection with $\norm{Q_n}\le C$.
  
  We are going to define sequences $(\widehat{E}_n)$,
  $(\widehat{F}_n)$ and $(\widehat{Q}_n)$ which satisfy all the
  properties described in the previous paragraph and, in addition,
  there exists a strictly increasing sequence $(m_n)$ in $\mathbb N$
  such that the following four conditions are satisfied
  \begin{enumerate}
    \item\label{leftEF} $m_{n-1}<\min\supp\widehat{E}_n$ and
                       $m_{n-1}<\min\supp\widehat{F}_n$;
    \item\label{leftQ} $\widehat{Q}_ny=0$ whenever $\supp y\le m_{n-1}$;
    \item\label{rightE} $m_n\ge\max\supp\widehat{E}_n$;
    \item\label{rightQ} $\norm{\widehat{Q}_ny}\le 2^{-n}\norm{y}$
                        whenever $\min\supp y>m_n$.
  \end{enumerate}
  We construct the sequences inductively. Let $m_0=0$, and suppose
  that we already constructed $\widehat{E}_i$, $\widehat{F}_i$,
  $\widehat{Q}_i$, and $m_i$ for all $i<n$.  Let $G$ and $G'$ be the
  subspaces of $\ell_p$ and $\ell_q$ respectively, consisting of all
  the vectors whose first $m_{n-1}$ coordinates are zero. Put
  $k=2m_{n-1}+n$. It follows from $\dim F_k=k$ and $\codim G'=m_{n-1}$ that
  $m_{n-1}+n\le\dim F_k\cap G'=\dim R^{-1}\bigl(F_k\cap G'\bigr)$ because
  $R_{|E_k}$ is an isomorphism. Since $\codim G=m_{n-1}$ we have
  $G\cap R^{-1}\bigl(F_k\cap G'\bigr)\ge n$. Let $\widehat{E}_n$ be an
  $n$-dimensional subspace of $G\cap R^{-1}\bigl(F_k\cap G'\bigr)$,
  and $\widehat{F}_n=R(\widehat{E}_n)$. Then $\widehat{E}_n\subseteq
  G$ and $\widehat{F}_n\subseteq G'$, hence (\ref{leftEF}) is
  satisfied.  Clearly, $\widehat{F}_n$ is $\widehat{C}$-complemented
  in $\ell_q$, where $\widehat{C}=C^2$.  Then there exists a
  projection $Q'\colon\ell_q\to\widehat{F}_n$ such that
  $\norm{Q'}\le\widehat{C}$. Let $\widehat{Q}_n=Q'P$, where $P$ is the
  basis projection of $\ell_q$ onto $[f_i]_{i\ge m_{n-1}}$.  Then
  $\widehat{Q}_n$ is again a projection from $\ell_q$ onto
  $\widehat{F}_n$, $\norm{\widehat{Q}_n}\le\widehat{C}$, and
  (\ref{leftQ}) is satisfied.  Since $\rank\widehat{Q}_n=n$, we can
  write $\widehat{Q}_n=\sum_{j=1}^nz_j\otimes z_j^*$, where
  $z_1,\dots,z_j\in\ell_p$ and $z_1^*,\dots,z_j^*\in\ell_q^*$. Then we
  can find $r\in\mathbb N$ sufficiently large, such that if
  $\norm{y}\le 1$ and $\min\supp y>r$ then $\abs{z_j^*(y)}$ is
  sufficiently small for all $j=1,\dots,n$, so that
  $\norm{\widehat{Q}y}\le 2^{-n}$. Let $m_n=\max\{r,s\}$, where
  $s=\max\supp\widehat{E}_n$, then (\ref{rightE}) and (\ref{rightQ})
  are satisfied.

  For convenience, we relabel $\widehat{E}_n$, $\widehat{F}_n$,
  $\widehat{Q}_n$, and $\widehat{C}$ as $E_n$, $F_n$, $Q_n$, and
  $C$ again.
  For every $n$ suppose that $V_n$ is a $C$-isomorphism of $\ell_2^n$
  onto $E_n$ with $\norm{V_n}=1$ and $\norm{V_n^{-1}}\le C$. Put
  $$V=\bigoplus_{n=1}^\infty V_n\colon 
    \Bigl(\bigoplus_{n=1}^\infty\ell_2^n\Bigr)_p\to 
  \Bigl(\bigoplus_{n=1}^\infty E_n\Bigr)_p.$$
  Since $E_n$'s are
  disjointly supported, we can consider
  $\bigl(\bigoplus_{n=1}^\infty E_n\bigr)_p$ as a subspace of
  $\ell_p$. It follows that $V$ is a $C$-isomorphism
  between $\bigl(\bigoplus_{n=1}^\infty\ell_2^n\bigr)_p$ and a
  subspace of $\ell_p$. Define 
  $$W\colon\ell_q\to\Bigl(\bigoplus_{n=1}^\infty\ell_2^n\Bigr)_q
    \text{ via }
    W\colon x\mapsto\Bigl(V_n^{-1}\bigl(R_{|E_n}\bigr)^{-1}
         Q_nx\Bigr)_{n=1}^\infty.$$
  We claim that $W$ is bounded. Indeed, pick $x\in\ell_q$. Then
  \begin{equation}\label{normW}
    \norm{Wx}=
    \Bigl(\sum_{n=1}^\infty\bignorm{V_n^{-1}\bigl(R_{|E_n}\bigr)^{-1}
             Q_nx}_2^q\Bigr)^{\frac{1}{q}}\le
    C^2\Bigl(\sum_{n=1}^\infty\norm{Q_nx}^q\Bigr)^{\frac{1}{q}}.
  \end{equation}
  Let $P_k$ be the basis projection from $\ell_q$ onto
  $[f_i]_{i=m_{k-1}+1}^{m_k}$. Then
  $x=\sum_{k=1}^\infty P_kx$. It follows from
  (\ref{leftQ}) that $Q_nP_kx=0$ whenever $k<n$. Furthermore,
  (\ref{rightQ}) yields 
  $\bignorm{Q_n\bigl(\sum_{k>n}P_kx\bigr)}\le 2^{-n}\norm{x}$. Also,
  $\norm{Q_nP_nx}\le C\norm{P_nx}$.
  Therefore, $\norm{Q_nx}\le C\norm{P_nx}+2^{-n}\norm{x}$.
  Using Cauchy-Schwartz inequality, we get
  \begin{displaymath}
   \Bigl(\sum_{n=1}^\infty\norm{Q_nx}^q\Bigr)^{\frac{1}{q}}\le
   \Bigl(\sum_{n=1}^\infty \bigl(C\norm{P_nx}\bigr)^q\Bigr)^{\frac{1}{q}}+
   \Bigl(\sum_{n=1}^\infty\bigl(2^{-n}\norm{x}\bigr)^q\Bigr)^{\frac{1}{q}}
   \le (C+1)\norm{x}.
  \end{displaymath}
  Together with (\ref{normW}) this yields that $W$ is bounded.

  Finally, it is easy to see that $WRV=\I2pq$, it follows easily that
  every \T{} operator factors through~$R$.
\end{proof}

\begin{corollary}
  Let $1<p\le 2\le q<\infty$. If $R\in\Lpq$ is not FSS, then
  $\J2\subseteq\iJ^R$.
\end{corollary}

\section{Operators not factorable through $\ell_2$}
\label{sec:walsh}

We employ the following known theorem (see
{\cite[Theorem~9.13]{Diestel:95}} or~\cite[Theorem~27.1]{Tomczak:89})
to deduce conditions for an operator in $\Lpq$ to factor
through~$\ell_r$.

\begin{theorem}\label{factor}
  Let $1\le r<\infty$, let $U\colon X\to Y$ be a bounded linear
  operator between the Banach spaces $X$ and~$Y$, and let $C\ge 0$.
  The following are equivalent:
  \begin{enumerate}
  \item\label{fact} There exists a subspace $L$ of $L_r(\mu)$, $\mu$ a
    measure, and a factorization $U=V\circ W$, where $V\colon L\to Y$
    and $W\colon X\to L$ are bounded linear operators with
    $\norm{V}\cdot\norm{W}\le C$.
  \item\label{domin} Whenever the finite sequences $(x_i)_{i=1}^n$ and
    $(z_i)_{i=1}^m$ in $X$ satisfy
    \begin{displaymath}
      \sum_{i=1}^m \bigabs{\langle x^*,z_i\rangle}^r\le 
      \sum_{i=1}^n \bigabs{\langle x^*,x_i\rangle}^r
          \mbox{ for all }x^*\in X^*,
      \mbox{then }
      \sum_{i=1}^m\norm{Uz_i}^r\le C^r \sum_{i=1}^n\norm{x_i}^r.
    \end{displaymath}
  \end{enumerate}
\end{theorem}

Let us use Theorem \ref{factor}  to state a criterion for an operator
$U:\ell_p^m\to\ell_q^m$ not to factor as $U=AB$ with
$\norm{B}_{p,r}\cdot\norm{A}_{r,q}\le C$.

\begin{corollary}\label{C:2}
  Let $m\in\mathbb N$, $C>1$, and $r>1$, and
  assume that
  $U$ is an invertible $m$ by $m$ matrix.
  Let $\delta=\norm{U^{-1}}_{r',r'}$.
  Then 
  $\norm{B}_{p,r}\cdot\norm{A}_{r,q}\ge \delta^{-1}$
  for any factorization $U=AB$.
  Moreover, if $\widetilde U$ is another $m$ by $m$ matrix with 
  \begin{equation}\label{approx}
    \norm{\widetilde U- U}_{p,q}\le
    \bigl(2\max\limits_{1\le i\le m}\norm{U^{-1}e_i}_p\bigr)^{-1},
  \end{equation}
  then it follows that for any factorization $\widetilde U=AB$ we have
  $\norm{B}_{p,r}\cdot\norm{A}_{r,q}\ge (2\delta)^{-1}$.
\end{corollary}

\begin{proof}
  For $i=1,\ldots,m$ we choose $x_i=e_i$ and $z_i=\delta^{-1} U^{-1} e_i$ and
  observe that for any $x^*\in\mathbb F^m$:
  \begin{multline*}
    \Bigl(\sum_{i=1}^m\bigabs{\langle x^*,z_i\rangle}^r\Bigr)^{1/r}
    =\delta^{-1}\Bigl(\sum_{i=1}^m\bigabs{\langle(U^{-1})^*x^*,e_i\rangle}^r\Bigr)^{1/r}
    =\delta^{-1}\norm{(U^{-1})^*x^*}_r \\
    \le\delta^{-1}\norm{U^{-1}}_{r',r'}\norm{x^*}_r
    =\Bigl(\sum_{i=1}^m\bigabs{\langle x^*,x_i\rangle}^r\Bigr)^{1/r},
  \end{multline*}
  which implies that the hypothesis of~(\ref{domin}) in
  Theorem~\ref{factor} is satisfied.  Secondly it follows that
  \begin{equation}\label{E:2.2}
    \sum_{i=1}^m\norm{Uz_i}^r_q=\delta^{-r}m=
        \delta^{-r}\sum_{i=1}^m\norm{x_i}^r_p,
  \end{equation}
  which means that the conclusion of~(\ref{domin}) in
  Theorem~\ref{factor} is not satisfied for any $C<\delta^{-1}$. It
  follows that condition~(\ref{fact}) in Theorem~\ref{factor} fails
  whenever $C<\delta^{-1}$.

  Now assume that $\widetilde U$ is another $m$ by $m$ matrix
  satisfying~(\ref{approx}), then it follows for $i=1,\ldots,m$ that
  \begin{eqnarray*}
    \norm{\widetilde U(z_i)}_q&\ge&\norm{U(z_i)}_q
       -\norm{(U-\widetilde U)(z_i)}_q\\
    &\ge&\tfrac{1}{2}\norm{U(z_i)}_q+
    \Bigl(\tfrac{1}{2}\norm{U(z_i)}_q
       -\norm{U-\widetilde U}_{p,q}\norm{z_i}_p\Bigr)\\
    &=&\tfrac{1}{2}\norm{U(z_i)}_q+\Bigl(\tfrac{1}{2\delta}-
    \norm{U-\widetilde U}_{p,q}\delta^{-1}\norm{U^{-1}e_i}_p\Bigr)\ge
        \tfrac{1}{2}\norm{U(z_i)}_q,
  \end{eqnarray*}
  which implies together with~(\ref{E:2.2}) that for $\widetilde U$ the
  conclusion of~(\ref{domin}) in Theorem \ref{factor} is not satisfied
  for any $C<\delta^{-1}/2$, hence~(\ref{fact}) fails in this case.
\end{proof}

\bigskip

We will now define an operator which will be cruical for the rest of the
paper. The following notations will be used throughout the rest of this
paper. Let $H_n$ be the $n$-th Hadamard matrix. That is, $H_1=(1)$,
$H_{n+1}=\left(
\begin{smallmatrix}
H_n & H_n\\
H_n & -H_n 
\end{smallmatrix}
\right)$ for every $n\ge 1$. Then $H_n$ is an $N\times N$ matrix where
$N=2^n$.
We view $\ell_p=\bigl(\bigoplus_{n=1}^\infty X_n\bigr)_p$ and
 $\ell_q=\bigl(\bigoplus_{n=1}^\infty Y_n\bigr)_q$, where
 $X_n=\ell_p^{2^n}$ and $Y_n=\ell_q^{2^n}$ are block subspaces of
 $\ell_p$ and $\ell_q$ respectively. We view $H_n$ as
 an operator from $X_n$ to $Y_n$. Define
\begin{equation}\label{U-def}
 U_n=N^{-\frac{1}{\min\{p',q\}}}H_n\mbox{ where }N=2^n,
 \mbox{ and let }U=\bigoplus_{n=1}^\infty U_n\colon\ell_p\to\ell_q.
\end{equation}

\begin{remark}\label{H-norms}
  Observe that $N^{-\frac{1}{2}}H_n$ is a unitary matrix on
  $\ell_2^N$. In particular, it is an isometry on $\ell_2^N$, hence
  $\norm{H_n}_{2,2}=N^{\frac{1}{2}}$, and $H_n^2=NI$.  One can easily
  verify that $\norm{H_n}_{1,\infty}=1$ and
  $\norm{H_n}_{1,1}=\norm{H}_{\infty,\infty}=N$.
\end{remark}

\begin{theorem}\label{U}
  If $p\le 2\le q$, then the operator $U$ defined by~(\ref{U-def}) has
  the following properties.
  \begin{enumerate}
    \item\label{norm-one} $\norm{U}_{p,q}=1$.
    \item\label{noncomp} $U$ is not compact.
    \item\label{U-fss} If $p'\neq q$ then $U$ is FSS.
    \item\label{r-fact} Let $p\le r\le q$. Then $U$ factors through
      $\ell_r$ when $p\le r\le q'$ or $p'\le r\le q$; otherwise $U\notin\Jr$.
    \item\label{notJ2} In particular, if $p\ne q$ then $U\notin\J2$.
  \end{enumerate}
\end{theorem}

\begin{proof}
  Using Riesz-Thorin Interpolation (e.g., \cite{Lindenstrauss:79})
  between $H_n$ as operator in $L(\ell_1,\ell_\infty)$ and as operator
  in $L(\ell_2,\ell_2)$, and using Remark~\ref{H-norms}, we obtain
  $\norm{H_n}_{r,r'}\le N^{\frac{1}{r'}}$ whenever $1\le r\le 2$.
  Similarly, interpolating between $\norm{H}_{1,1}$ and
  $\norm{H_n}_{2,2}$, and between $\norm{H_n}_{2,2}$ and
  $\norm{H}_{\infty,\infty}$ we obtain $\norm{H}_{r,r}\le
  N^{\frac{1}{\min\{r,r'\}}}$ whenever $1\le r\le\infty$.

  Define $U^{(r)}_n=N^{-\frac{1}{r'}}H_n$ and
  $U^{(r)}=\bigoplus_{n=1}^\infty U_n^{(r)}$, then
  $\norm{U_n^{(r)}}_{r,r'}\le 1$ for every~$n$, hence
  $\norm{U^{(r)}}_{r,r'}\le 1$.
  Viewing $U$ as an operator in $\Lpq$, we can write
  \begin{equation}
    \label{U-cases}
    U=
    \begin{cases}
      \ell_p\xrightarrow{U^{(p)}}\ell_q & \mbox{ when }p'=q,\\
        \ell_p\xrightarrow{U^{(p)}}\ell_{p'}\xrightarrow{I_{p',q}}\ell_q
      & \mbox{ when }p'<q,\mbox{ and}\\
      \ell_p\xrightarrow{I_{p,q'}}\ell_{q'}\xrightarrow{U^{(q')}}\ell_q
      & \mbox{ when }p<q'.
    \end{cases}
  \end{equation}
  It follows immediately that $\norm{U}_{p,q}\le 1$.
  Since $\fss$ is an ideal, (\ref{U-fss})~follows from
  Proposition~\ref{Ipq-fss}. It also follows from~(\ref{U-cases}) that
  $U$ factors through $\ell_r$ if $p\le r\le q'$ or $p'\le r\le q$,

  Consider first the case $p'\le q$. Then $U_n=N^{-\frac{1}{p'}}H_n$.
  Let $h_{n,i}=H_ne_i$, the $i$-th column of
  the $n$-th Hadamard matrix.
  It follows from $H_n^2=NI$ that
  $U_nh_{n,i}=N^{-\frac{1}{p'}}H_n^2e_i=N^{\frac{1}{p}}e_i$.
  Thus, $\norm{U_nh_{n,i}}_q=N^{\frac{1}{p}}=\norm{h_{n,i}}_p$, so
  that $\norm{U_n}_{p,q}=1$. Hence, $U$ is not compact, and
  $\norm{U}_{p,q}=1$ by Remark~\ref{block-diag}.
  
  Next, suppose that $p<r<p'\le q$. We use Corollary~\ref{C:2} to show
  that $U\notin\Jr$ in this case. Indeed, assume to the contrary that
  $U\in\Jr$. Then there exists $\widetilde{U}$ such that
  $\norm{U-\widetilde{U}}<\frac{1}{2}$ and $\widetilde{U}$ factors
  through~$\ell_r$. Let $C$ be the $\ell_r$\hyph factorization
  constant of~$\widetilde{U}$. Since $p<\min\{r,r'\}$ one can choose
  $n$ so that $C<\tfrac{1}{2}N^{\frac{1}{p}-\frac{1}{\min\{r,r'\}}}$,
  where $N=2^n$. Let $\widetilde{U}_n$ be the $N\times N$ submatrix of
  $\widetilde{U}$ corresponding to the $n$-th block of~$U$, that is,
  $\widetilde{U}_n=Q_n\widetilde{U}P_n$, where $P_n$ (respectively,
  $Q_n$) is the canonical projection from $\ell_p$ (respectively,
  $\ell_q$) onto the span of $e_{N+1},\dots,e_{2N}$. Then the
  $\ell_r$\hyph factorization constant of $\widetilde{U}_n$ is at
  most~$C$.  It follows from
  $\norm{U_n^{-1}e_i}_p=\norm{N^{-\frac{1}{p}}h_{n,i}}_p=1$ that
  $$\bignorm{U_n-\widetilde{U}_n}\le\bignorm{U-\widetilde{U}}<\tfrac{1}{2}=
    \bigl(2\max\limits_{1\le i\le N}\bignorm{U_n^{-1}e_i}_p\bigr)^{-1}.$$
  Let $\delta=\norm{U_n^{-1}}_{r',r'}$. It follows from $H_n^2=NI$ and
  $U_n=N^{-\frac{1}{p'}}H_n$ that $U_n^{-1}=N^{-\frac{1}{p}}H$, so that
  $$\delta=N^{-\frac{1}{p'}}\norm{H_n}_{r',r'}\le
    N^{-\frac{1}{p}+\frac{1}{\min\{r,r'\}}}.$$
  Corollary~\ref{C:2} yields that the factorization constant of
  $\widetilde{U}_n$ is at least
  $(2\delta)^{-1}\ge\frac{1}{2}N^{\frac{1}{p}-\frac{1}{\min\{r,r'\}}}>C$,
  contradiction.
  
  The case $p<q'$ can be reduced to the previous case by duality.
  Indeed, it follows from~(\ref{U-cases}) that
  $U^*=I_{q,p'}U^{(q')}\colon\ell_{q'}\to\ell_{p'}$. It follows that
  if $p\le r\le q'$ then $I_{q,p'}$ and, therefore, $U^*$ factors
  through $\ell_{r'}$. Hence, $U$ factors through~$\ell_r$.
  Furthermore, since $H_n$ is symmetric for every $n$, it follows that
  $U_n^*$ coincides with $U_n$ as a matrix and
  $\norm{U_n^*}_{q',p'}=1$.  Applying the previous argument, we
  observe that $U^*$ is non-compact and $\norm{U^*}_{q',p'}=1$, hence
  the same is true for $U$.  Furthermore, if $q'<r<q$, then
  $U^*\notin\iJ^{\ell_{r'}}$ so that $U\notin\Jr$.

  Finally, (\ref{notJ2}) follows immediately from~(\ref{r-fact}).
\end{proof}

\begin{remark}
  If $p<r<r'<q$ then the operator $\widetilde U$ defined as
  $$\ell_p\xrightarrow{I_{p,r}}\ell_r\xrightarrow{U^{(r)}}
    \ell_{r'}\xrightarrow{I_{r',q}}\ell_q$$
  is compact. Indeed, as a matrix
  $$\widetilde{U}_n=U_n^{(r)}=
     N^{-\frac{1}{r'}}H_n=N^{\frac{1}{\min\{p',q\}}-\frac{1}{r'}}U_n.$$
  It follows from $\norm{U_n}_{p,q}=1$ and $r'<\min\{p',q\}$ that
  $\norm{\widetilde U_n}_{p,q}=N^{\frac{1}{\min\{p',q\}}-\frac{1}{r'}}\to 0$
  as $n\to 0$.
\end{remark}

\begin{remark}
  It follows from Theorem~\ref{U}(\ref{r-fact}) that $\Jr$ is proper
  when $\max\{p,q'\}<r<\min\{p',q\}$. In particular, $\J2$ is proper.
  It follows from Remark~\ref{Iprq} and Theorem~\ref{U}(\ref{r-fact})
  that $\J2\subsetneq\Jr$ whenever $p<r<q'$ or $p'<r<q$.  We do not
  know, however, whether $\Jr$ is proper in this case.
\end{remark}

\medskip

Next, we show that if $U'$ is another ``$U$-like'' operator then $U$
and $U'$ factor through each other.

Again, we view $\ell_p=\bigl(\bigoplus_{n=1}^\infty X_n\bigr)_p$ and
$\ell_q=\bigl(\bigoplus_{n=1}^\infty Y_n\bigr)_q$, where
$X_n=\ell_p^{2^n}$ and $Y_n=\ell_q^{2^n}$. Denote the basis vectors in
$X_n$ and $Y_n$ by $e^{(n)}_1,\dots,e^{(n)}_{2^n}$ and
$f^{(n)}_1,\dots,f^{(n)}_{2^n}$, respectively. We can view $H_n$ and
$U_n$ as operators from $X_n$ to~$Y_n$.

\begin{theorem}
  Suppose that $(n_i)$ is an increasing sequence, and let
  $\widetilde{U}=\bigoplus_{i=1}^\infty U_{n_i}$, viewed as an
  operator from $\ell_p=\bigl(\bigoplus_{i=1}^\infty X_{n_i}\bigr)_p$
  to $\ell_q=\bigl(\bigoplus_{i=1}^\infty Y_{n_i}\bigr)_q$.  Then $U$
  and $\widetilde{U}$ factor through each other.
\end{theorem}

\begin{proof}
  Consider the following diagram
  $$\ell_p=\bigl(\toplus_{i=1}^\infty X_{n_i}\bigr)_p
  \overset{\imath}{\hookrightarrow}
  \bigl(\toplus_{n=1}^\infty X_n\bigr)_p\xrightarrow{U}
  \bigl(\toplus_{n=1}^\infty Y_n\bigr)_q\xrightarrow{R}
  \bigl(\toplus_{i=1}^\infty Y_{n_i}\bigr)_q=\ell_q,$$
  where $\imath$ is the canonical embedding, and $R$ is the canonical
  projection. We can view $\imath$ and  $R$ as operators on $\ell_p$
  and $\ell_q$ respectively. Thus, we get $\widetilde{U}=RU\imath$.

  Next, we prove that $U$ factors through~$\widetilde{U}$. First, we
  show that whenever $n<m$ then there exists operators $C\colon X_n\to X_m$
  and $D\colon Y_m\to Y_n$ such that $U_n=DU_mC$ and
  $\norm{C}_{p,p}\le 1$ and $\norm{D}_{q,q}\le 1$.

  First, we consider the case $q\le p'$.
  Define $C_n\colon X_n\to X_{n+1}$ via $C_ne^{(n)}_i=e^{(n+1)}_i$ as
  $i=1,\dots,2^n$. Clearly, $C_n$ is an isometry.
  
  Let $Z_n$ be the subspace of $Y_{n+1}$ consisting of all the vectors
  whose first half coordinates are equal to the last half coordinates
  respectively, that is, $Z_n=\Span\{f^{(n+1)}_i+f^{(n+1)}_{i+2^n}\mid
  i=1,\dots,2^n\}$.  Let $P_n$ be the ``averaging'' projection from
  $Y_{n+1}$ onto $Z_n$ given by
  $$P_n\Bigl(\sum\limits_{i=1}^{2^{n+1}}\alpha_if^{(n+1)}_i\Bigr)=
    \sum\limits_{i=1}^{2^n}\frac{\alpha_i+\alpha_{i+2^n}}{2}
    \bigl(f^{(n+1)}_i+f^{(n+1)}_{i+2^n}\bigr).$$
  Then $\norm{P_n}=1$.
  
  Define $B_n\colon Z_n\to Y_n$ via
  $B(f^{(n+1)}_i+f^{(n+1)}_{i+2^n})=2^{\frac{1}{q}}f^{(n)}_i$, then
  $B_n$ is an isometry. Hence, $D_n=B_nP_n\colon Y_{n+1}\to Y_n$ is of
  norm one.
  
  Fix $1\le i\le 2^n$. Since $C_ne^{(n)}_i=e^{(n+1)}_i$, 
  $H_{n+1}C_ne^{(n)}_i$ is the $i$-th column of
  $H_{n+1}$. Since $i\le 2^n$ it follows from the construction of
  $H_n$'s that the $i$-th column of $H_{n+1}$ is exactly the $i$-th
  column of $H_n$ repeated twice. In particular,
  $H_{n+1}C_ne^{(n)}_i\in Z_n$ and, therefore,
  $H_{n+1}C_ne^{(n)}_i=P_nH_{n+1}C_ne^{(n)}_i$. Finally,
  $$B_nP_nH_{n+1}C_ne^{(n)}_i=2^{\frac{1}{q}}(\mbox{the $i$-th column
    of $H_n$})=2^{\frac{1}{q}}H_ne^{(n)}_i.$$
  It follows that
  $D_nH_{n+1}C_n=2^{\frac{1}{q}}H_n$. It follows from
  $H_n=2^{\frac{n}{q}}U_n$ that $D_nU_{n+1}C_n=U_n$. Iterating this
  $m-n$ times, we get $DU_mC=U_n$ where $C\colon X_n\to X_m$ is an
  isometry, $D\colon Y_m\to Y_n$ is of norm one.

  If $q\ge p'$, then we consider the adjoint operators. Note that 
  $U_n^*=U_n$ as matrices. Applying the
  previous argument we find matrices $C$ and $D$ such that
  $U^*_n=DU_m^*C$ with $\norm{C}_{q',q'}\le 1$ and
  $\norm{D}_{p',p'}\le 1$. Then $U_n=C^*U_mD^*$ is a required
  factorization in the case $q\ge p'$.

  It follows that for every $i$ we have
  \begin{equation}
    \label{eq:tildes}
    \widetilde D_iU_{n_i}\widetilde C_i=U_i
  \end{equation}
  for some contractions $\widetilde C_i\colon X_i\to X_{n_i}$ and
  $\widetilde D_i\colon X_{n_i}\to X_i$. Let
  $$\widetilde C=
    \bigl(\oplus_{i=1}^\infty\widetilde C_i\bigr)
    \colon\bigl(\oplus_{i=1}^\infty X_i\bigr)_p
    \to\bigl(\oplus_{i=1}^\infty X_{n_i}\bigr)_p$$
  and
  $$\widetilde D= \bigl(\oplus_{i=1}^\infty\widetilde D_i\bigr)\colon
  \bigl(\oplus_{i=1}^\infty X_{n_i}\bigr)_q
  \to\bigl(\oplus_{i=1}^\infty X_i\bigr)_q.$$
  Then $\widetilde
  C\colon\ell_p\to\ell_p$ and $\widetilde D\colon\ell_q\to\ell_q$ are
  bounded, and by (\ref{eq:tildes}) we have $\widetilde
  D\widetilde{U}\widetilde C=U$.
\end{proof}

It follows that any two operators of type $\widetilde{U}$ generated by different
sequences factor through each other.

\section{The operator $U$ is FSS}\label{Ufss}

Again, let $U$ be the operator defined by~(\ref{U-def}).
Theorem~\ref{U}(\ref{U-fss}) states that $U$ is FSS when $p\neq q'$.
We will show in this section that $U$ is still FSS when $1<p=q'$. The
argument requires certain preparation.

Recall that the $n$-th $s$-number $s_n(T)$ of an operator $T\in L(H)$
on a Hilbert space $H$ is defined as the distance from $T$ to the set
of all operators in $L(H)$ of rank $n-1$. For $1\le r<\infty$, the
Schatten norm $\bignorm{T}_{S_r}$ of $T$ equals the $\ell_r$ norm of
the sequence of the $s$-numbers. We say that $T$ belongs to Schatten
class $S_r$ if $\bignorm{T}_{S_r}<\infty$. We denote by $S_\infty$ the
set of all compact operators equipped with the operator norm.

\begin{lemma}\label{Sh-FSS}
  If $T\in L(H)$ such that $\bignorm{T}_{S_q}=1$ and $\inf\limits_{x\in
    F,\norm{x}=1} \norm{Tx}\ge\varepsilon$ for a subspace $F$ of~$H$,
  then $\dim F\le\varepsilon^{-q}$.
\end{lemma}

\begin{proof}
  Suppose that $\dim F=k$. For every operator $S$ of rank $k-1$ there
  exists $x\in F$ such that $\norm{x}=1$ and $Sx=0$. It follows that
  $\norm{T-S}\ge\norm{Tx}\ge\varepsilon$, so that
  $s_1\ge\dots\ge s_k\ge\varepsilon$. Therefore,
  $1=\bignorm{T}_{S_q}^q\ge k\varepsilon^q$. Hence
  $k\le\varepsilon^{-q}$.
\end{proof}

We will also utilize the following result of
Maurey,~\cite[Corollary~11, p.~21]{Maurey:74}.

\begin{theorem}\label{Maurey}
  Let $(\Omega,\mu)$ be a measure space, $X$ and $Y$ two quasi-normed
  vector spaces, $0<u\le v<\infty$,
  $\frac{1}{u}=\frac{1}{v}+\frac{1}{r}$, 
  $T$ a bounded operator from a closed subspace $E$ of
  $L_v(\Omega;X)$ to~$Y$, and $C>0$. Then the following
  are equivalent.
  \begin{enumerate}
   \item\label{M-fact}
     There exists a closed subspace $F$ of $L_u(\Omega;X)$ such that
     $T$ factors $T=V\circ M_g$, where $V\colon F\to Y$ with
     $\norm{V}\le C$, and $M_g$ a multiplication operator with
     multiplier $g\in L_r(\mu)$ with $\norm{g}_r\le 1$.
   \item\label{M-sums}
     For any $x_1,\dots,x_n$ in $E$,
     \begin{displaymath}
      \Bigl(\sum\limits_{i=1}^n\norm{Tx_i}^u\Bigr)^{\frac{1}{u}}\le
      C\biggl[\int\Bigl(\sum\limits_{i=1}^n\norm{x_i}_X^u\Bigr)
      ^{\frac{v}{u}}\,d\mu\biggr]^{\frac{1}{v}}.
     \end{displaymath}
  \end{enumerate}
\end{theorem}

\begin{corollary}\label{Maurey-rev}
  Let $(\Omega,\mu)$ be a measure space. Suppose that $q=p'$ and
  $\frac{1}{p}=\frac{1}{2}+\frac{1}{r}$.
  \begin{enumerate}
  \item\label{M:q2} If $T\colon L_q(\Omega)\to\ell_2^k$ then $T$ can
    be factored through a multiplication operator on $L_2(\Omega)$ as
    $T=SM_g$, where $S\colon L_2(\Omega)\to\ell_2^k$ with $\norm{S}\le
    K_G\norm{T}$ and $\norm{g}_r=1$.
  \item\label{M:2p} If $T\colon\ell_2^k\to L_p(\Omega)$ then $T$ can
    be factored through a multiplication operator on~$L_2$, that is,
    $T=M_hS$, where $S\colon\ell_2^k\to L_2(\Omega)$ with $\norm{S}\le
    K_G\norm{T}$ and $\norm{h}_r\le 1$.
  \end{enumerate}
\end{corollary}

\begin{proof}
  Suppose that $T\colon L_q(\Omega)\to\ell_2^k$.  We verify that
  condition~(\ref{M-sums}) of Theorem~\ref{Maurey} holds for
  $X=\mathbb F$, $u=2$, $v=q=p'$, and $r>1$ such that
  $\frac{1}{p}=\frac{1}{2}+\frac{1}{r}$ (which is equivalent to
  $\frac{1}{2}=\frac{1}{v}+\frac{1}{r}$). Let $f_1,\dots,f_n\in L_q$.
  Then 
  \begin{displaymath}
    \sum\limits_{i=1}^n\bignorm{Tf_i}^2=
    \sum\limits_{i=1}^n\sum\limits_{j=1}^k\bigabs{(Tf_i)_j}^2=
    \sum\limits_{j=1}^k\sum\limits_{i=1}^n\bigabs{(Tf_i)_j}^2=
    \Bignorm{\Bigl(\sum\limits_{i=1}^n\bigabs{Tf_i}^2\Bigr)
       ^{\frac{1}{2}}}_{\ell_2}^2,
  \end{displaymath}
  where the last expression is the norm of the sequence
  $\Bigl(\bigl(\sum_{i=1}^n\abs{(Tf_i)_j}^2\bigr)^
   {\frac{1}{2}}\Bigr)_{j=1}^n$.
  It follows from~\cite[Theorem~1.f.14]{Lindenstrauss:79} that
  \begin{displaymath}
    \Bignorm{\Bigl(\sum\limits_{i=1}^n\bigabs{Tf_i}^2\Bigr)
       ^{\frac{1}{2}}}_{\ell_2}\le
    K_G\norm{T}\Bignorm{\Bigl(\sum\limits_{i=1}^n\abs{f_i}^2\Bigr)
       ^{\frac{1}{2}}}_{L_q}=
    K_G\norm{T}\biggl[\int\Bigl(\sum\limits_{i=1}^n\abs{f_i}^2\Bigr)
       ^{\frac{q}{2}}\,d\mu\biggr]^{\frac{1}{q}},
  \end{displaymath}
  where $K_G$ is the Grothendieck constant.
  Now~(\ref{M:q2}) follows from Theorem~~\ref{Maurey}.
  To prove~(\ref{M:2p}), apply (\ref{M:q2}) to~$T^*$.
\end{proof}

\begin{lemma}\label{HS}
  Consider a product of three bounded operators
  \begin{displaymath}
    S\colon L_2\xrightarrow{M_{\psi}}L_1\xrightarrow{T}\ell_\infty
    \xrightarrow{D}\ell_2,
  \end{displaymath}
  where $\psi\in L_2$ and $D=\diag(d_j)_{j=1}^\infty$.
  Then the Hilbert-Schmidt norm
  $\norm{S}_{\mathrm{HS}}\le\norm{\psi}_2\norm{T}\bignorm{(d_j)}_2$.
\end{lemma}

\begin{proof}
  Observe that
  \begin{displaymath}
    S\colon f\mapsto\psi f\mapsto
    \bigl(\langle g_n,\psi f\rangle\bigr)_{n=1}^\infty\mapsto
    \bigl(d_n\langle g_n,\psi f\rangle\bigr)_{n=1}^\infty
  \end{displaymath}
  for some uniformly bounded sequence $(g_n)_{n=1}^\infty$ in
  $L_\infty$, and $\norm{T}=\sup_n\norm{g_n}_\infty$. Let
  $(f_i)_{i=1}^\infty$ be an orthonormal basis of~$L_2$, then
  \begin{multline*}
    \norm{S}_{\rm HS}^2=\sum\limits_{i=1}^\infty\norm{Sf_i}^2=
    \sum\limits_{i=1}^\infty\sum\limits_{n=1}^\infty
      d_n^2\langle g_n,\psi f_i\rangle^2=
    \sum\limits_{n=1}^\infty d_n^2
    \sum\limits_{i=1}^\infty\langle\psi g_n,f_i\rangle^2=\\
    \sum\limits_{n=1}^\infty d_n^2\cdot\norm{\psi g_n}_{L_2}^2\le
    \norm{\psi}_2\cdot\bigl(\sup\limits_n\norm{g_n}_\infty^2\bigr)
      \cdot\bignorm{(d_j)}_2^2.
  \end{multline*}
\end{proof}

By $L_p^N$ we denote the discrete $L_p$-space on a set of $N$ elements
endowed with the uniform probability measure. Clearly, $L_p^N$ is
isometrically isomorphic to $\ell_p^N$. We will view $L_p^N$ as a
subspace of $L_p=L_p[0,1]$ under the natural isometric embedding.
Namely, the \mbox{$i$-th} basis vector $e_i$ of $\ell_p^N$ corresponds to the
function $N^{\frac{1}{p}}\chi_{[\frac{i-1}{N},\frac{i}{N}]}$. In
particular, every $N\times N$ matrix can be viewed as an operator on
$L_p^N$.

\begin{theorem}[{\cite{Pis}}]\label{dimension}
  Suppose that $T\colon L_p^N\to\ell_q^N$ for some $1\le p<2$ and
  $q=p'$. Let $E$ be a $k$\hyph dimensional subspace of $L_p^N$, and
  $C_1$, $C_2$, and $C_3$ be positive constants such that
  \begin{enumerate}
  \item\label{i:norms} $\norm{T}_{L_2^N,\ell_2^N}\le 1$ and
        $\norm{T}_{L_1^N,\ell_\infty^N}\le 1$;
  \item $E$ is $C_1$-isomorphic to $\ell_2^k$;
  \item $F=T(E)$ is $C_2$-complemented in $\ell_q^N$; and
  \item $T_{|E}$ is invertible and $\bignorm{(T_{|E})^{-1}}\le C_3$ 
  \end{enumerate}
  then $k\le\bigl(C_1^3C_2C_3^2K_G^2\bigr)^q$.
\end{theorem}

\begin{proof}
  Suppose that $T$, $E$, and $F$ satisfy the hypotheses for some
  $C_1$, $C_2$, and~$C_2$. Let $r$ be such that
  $\frac{1}{p}=\frac{1}{2}+\frac{1}{r}$. There exists an
  isomorphism $V\colon\ell_2^k\to E$ such that $\norm{V}\le 1$ and
  $\norm{V^{-1}}\le C_1$. By Corollary \ref{Maurey-rev}(\ref{M:2p})
  $V$ factors through $L_2^N$. Namely, $V=M_gS$ such that
  $S\colon\ell_2^k\to L_2^N$ with $\norm{S}\le C_1K_G$ and
  $\norm{g}_r\le 1$. Let $J\colon E\to L_p^N$ be the canonical
  inclusion map. 
  \begin{displaymath}
    \begin{CD}
      L_2^N @ >M_g>> L_p^N @ >T>> \ell_q^N @ >D>{\mathrm{diagonal}}> \ell_2^N\\
      @ ASAA @ AJA{\mathrm{incl.}}A @ V{\mathrm{proj.}}VQV @ VVRV \\
      \ell_2^k @ >V>{C_1-\mathrm{isom.}}> E @ >{T_{|E}}>> F @ 
         >W>{C_1C_3-\mathrm{isom}}> \ell_2^k
    \end{CD}
  \end{displaymath}
  Let $Q$ be a projection from $\ell_q^N$ onto $F$ with $\norm{Q}\le
  C_2$. It follows from~(\ref{i:norms}) that
  $\norm{T}_{L_p^N,\ell_q^N}\le 1$. Then $F$ is
  $C_1C_3$-isomorphic to $\ell_2^k$. Let $W\colon
  F\to\ell_2^k$ be an isomorphism such that $\norm{W}\le 1$ and
  $\norm{W^{-1}}\le C_1C_3$. Corollary \ref{Maurey-rev}(\ref{M:q2})
  implies that $WQ$ factors through $\ell_2^N$, that is, $WQ=RD$ where
  $R\colon\ell_2^N\to\ell_2^k$ with $\norm{R}\le K_G\norm{WQ}\le
  C_2K_G$, and $D$ is a
  multiplication (or diagonal) operator $D=\diag(d_j)_{j=1}^N$ with
  $\bignorm{(d_j)}_{\ell_r^N}\le 1$.
  
  We proceed with a version of the classical complex interpolation
  argument (see, e.g., \cite{Bergh:76}). Let $Z=\{z\in\mathbb C\mid
  0\le\re z\le 1\}$, and define a function $F$ from $Z$ to the unit
  ball $B(L_2^N,\ell_2^N)$ of $L(L_2^N,\ell_2^N)$ as follows:
  \begin{equation}\label{function}
    F(z)=\abs{D}^{(1-z)\frac{r}{2}}\sign D T
         M_{\abs{g}^{(1-z)\frac{r}{2}}\sign g}.
  \end{equation}
  Here, as usually, $\abs{D}=\diag(\abs{d_j})$ and $\sign
  D=\diag(\sign d_j)$. Observe that $F$ is analytic in the interior of
  $Z$ as a function from $Z$ to $\mathbb C^N\times\mathbb C^N$.
  Furthermore, $F$ is continuous and bounded on~$Z$. A direct
  calculation shows that if $\frac{1}{r}=\frac{1-\theta}{2}$ then
  $F(\theta)=DTM_g$.
 
  If $\re z=1$, it follows from~(\ref{function}) that
  $F(1+it)=A_tTB_t$,
  where
  $A_t=\abs{D}^{-\frac{itr}{2}}\sign D$
  and
  $B_t= M_{\abs{g}^{-\frac{itr}{2}}\sign g}$.
  Notice that $A_t$ and $B_t$  viewed
  as operators from $\ell_2^N$ to $\ell_2^N$ and from $L_2^N$ to
  $L_2^N$ respectively are contractions. It follows that
  \begin{equation}\label{eq:Re1}
    \bignorm{F(z)}_{L_2^N,\ell_2^N}\le\norm{T}_{L_2^N,\ell_2^N}\le
     1\text{ whenever }\re z=1.
  \end{equation}
 
  If $\re z=0$ then we can write
  $$F(it)=A_t\abs{D}^{\frac{r}{2}}T
  M_{\abs{g}^{\frac{r}{2}}}B_t.$$ It can be easily verified that
  $\bignorm{\abs{g}^{\frac{r}{2}}}_{L_2^N}\le 1$ and
  $\bignorm{\bigl(\abs{d_i}^{\frac{r}{2}}\bigr)}_{\ell_2^N}\le 1$.
  Since $\norm{T}_{L_1^N,\ell_\infty^N}\le 1$, it follows by
  Lemma~\ref{HS} that
  \begin{equation}\label{eq:Re0}
    \bignorm{F(z)}_{\mathrm{HS}}\le 1\text{ whenever }\re z=0.
  \end{equation}
  Denote $S_q^N=S_q(L_2^N,\ell_2^N)$. It is known (see,
  e.g.~\cite[Theorem~13.1]{Krein:65}) that the Schatten classes
  interpolate like $L_p$-spaces. Since
  $$\tfrac{1}{\infty}(1-\theta)+ 
     \tfrac{1}{2}\theta=\tfrac{1}{2}-\tfrac{1}{r}=\tfrac{1}{q},$$
  it follows that $(S^N_\infty,S^N_2)_\theta= S^N_q$.
                                                                                
  On the other hand, by definition of a complex interpolation space,
  \begin{multline*}
    B_{(S^N_\infty,S^N_2)_{\theta}}=\Bigl\{f(\theta)\mid
    f\colon Z\to B(L_2^N,\ell_2^N)\mbox{ analytic,}\\
    \bignorm{f_{|\{\re z=0\}}}_{S_2}\le 1\mbox{ and }
    \bignorm{f_{|\{\re z=1\}}}_{S_\infty}\le 1\Bigr\}.
  \end{multline*}
  Since $\norm{\cdot}_{S_2}=\norm{\cdot}_{\mathrm{HS}}$
  and $\norm{\cdot}_{S_\infty}=\norm{\cdot}_{L_2^N,\ell_2^N}$, it follows
  from~(\ref{eq:Re1}) and~(\ref{eq:Re0}) that
  $DTM_g=F(\theta)\in B_{(S^N_\infty,S^N_2)_{\theta}}$ and, thus,
  $\norm{DTM_g}_{S_q^N}\le 1$.
  It follows that
  \begin{displaymath}
    \norm{WTV}_{S_q}=
    \bignorm{RDTM_gS}_{S_q}\le\norm{R}\bignorm{DTM_g}_{S_q}\norm{S}
    \le C_1C_2K_G^2.
  \end{displaymath}
  Note that $\bignorm{(WTV)^{-1}}\le C_1^2C_3^2$.
  It follows from Lemma~\ref{Sh-FSS} that
  \begin{displaymath}
    k\le
    \Bigl(\frac{1}{C_1^2C_3^2}/\bigl(C_1C_2K_G^2\bigr)\Bigr)^{-q}
    =\bigl(C_1^3C_2C_3^2K_G^2\bigr)^q.
  \end{displaymath}
\end{proof}

We also need the following lemma, which generalizes
Lemma~\ref{flat-max}. Assume that $X$ is a Banach space with an FDD
$(X_n)_{n=1}^\infty$. Let $P_n$ be the canonical projection from $X$
onto~$X_n$, and assume that $X$ satisfies the following condition,
which means that $X$ is far apart from a $c_0$-sum of the $X_n$'s:
\begin{equation}\label{non-c-0}
  \parbox{11cm}{\it for any $\delta>0$ there is a $k=k(\delta)$ in
  $\mathbb N$ so that whenever $x\in S_X$, then
  $\card\{ n\in\mathbb N: \norm{P_nx}\ge \delta\}<k$.}
\end{equation}
Suppose that for every $n\in\mathbb N$ we are given a seminorm $q_n$
on $X_n$ such that $q_n(x)\le\norm{x}$, where
$q_n(x)$ stands for $q_n(P_nx)$ whenever $x\in X$.

\begin{lemma}\label{flat-FDD}
  Suppose that $X$, $(X_n)$, and $(q_n)$ are as in the preceding
  paragraph and $0<r\le 1$. Then there exists $\varepsilon>0$ such
  that for every $l\in\mathbb N$ there exists $L\in\mathbb N$ such
  that for every $L$\hyph dimensional subspace $G$ of $X$ such that
  $\max_{n\in\mathbb N}q_n(x)\ge r\norm{x}$ for all $x\in G$ there
  exists an $l$-dimensional subspace $F\subseteq G$ and an index $n_0$
  such that $q_{n_0}(x)\ge\varepsilon\norm{x}$ for all $x\in F$.
\end{lemma}

To prove Lemma~\ref{flat-FDD} we need the following stabilization
result, see, e.g.,~\cite[p.6]{Milman:86}.

\begin{theorem}\label{stabilization}
  For every $n\in\mathbb N$, $\varepsilon>0$ and $c>0$ there is an
  $N=N(n,\varepsilon,c)\in\mathbb N$ so that for any $N$-dimsensional
  space $E$, and any Lipschitz map $f:S_E\to \mathbb R$ whose
  Lipschitz constant does not exceed $c$, there is an $n$\hyph dimensional
  subspace $F$ of $E$ so that
  $$\max \{f(x): x\in S_F\}- \min \{f(x): x\in S_F\}\le \varepsilon.$$
\end{theorem}
                  
\begin{proof}[Proof of Lemma \ref{flat-FDD}]
  Let $k(\cdot)$ be the function defined in~(\ref{non-c-0}). Put
  $$m=k\bigl(\tfrac{r^2}{4}\bigr),\quad\delta=\tfrac{r}{4m},\quad
  \mbox{ and }\quad s=k(\delta).$$
  It suffices to show that for $l'\in\mathbb N$ there exists $L$ so
  that if $G$ is a subspace of $X$ of dimension $L$ and
  $\max_{n\in\mathbb N}q_n(x)\ge r\norm{x}$ for all $x\in G$ then $G$
  has an $l'$\hyph dimensional subspace $F'$ and a set $I\subset
  \mathbb N$ with $\card I=s$ such that $\max_{n\in
    I}q_n(x)\ge\delta\norm{x}$ for all $x\in F'$.
  
  Indeed, once we prove this formaly weaker claim, we can take
  a number $l'$ large enough, so that we can apply
  Theorem~\ref{stabilization} $s$ times to deduce that $F'$ has an
  $l$\hyph dimensional subspace $F$, which has the property that for all
  $n\in I$
  $$\max_{x\in S_F} q_n(x)-\min_{x\in S_F} q_n(x) \le
  \tfrac{\delta}{2}.$$
  Now pick any $y\in S_F$, then $q_{n_0}(y)=\max_{n\in I}q_n(y)\ge\delta$
  for some $n_0\in I$. Then for every $x\in S_F$ we have
  $$q_{n_0}(x)\ge
    \min\limits_{z\in S_F}q_{n_0}(z)\ge
    \max\limits_{z\in S_F}q_{n_0}(z)-\tfrac{\delta}{2}\ge
    q_{n_0}(y)-\tfrac{\delta}{2}\ge\tfrac{\delta}{2},$$
  so that the statement of our Lemma is satisfied for
  $\varepsilon=\frac{\delta}{2}$.

  Let $l'\in\mathbb N$ and define numbers $L_0,L_1,\dots,L_m$ as
  follows. Put $L_0=l'$, and, assuming that $L_0,L_1,\dots,L_n$,
  $n<m$, have already been defined, we use Theorem~\ref{stabilization}
  to choose $L_{n+1}$ large enough so that for every $L_{n+1}$\hyph
  dimensional subspace $G$ of $X$ and every Lipschitz-1 map $f\colon
  S_G\to\mathbb R$ there is an $L_n$\hyph dimensional subspace
  $G'\subseteq G$ such that
  $$\max\limits_{x\in G'}f(x)-\min\limits_{x\in G'}f(x)\le\delta.$$
  
  Let $L=L_m$. Assume that out claim is false. This would mean that
  there exists a subspace $G$ of $X$ of $\dim G=L$ such that
  \begin{equation}
   \label{qn_ge_r}
    \max\limits_{n\in\mathbb N}q_n(x)\ge r\norm{x}\mbox{ for all }x\in
    G,\mbox{ and}
  \end{equation}
  \begin{equation}
   \label{s-set}
     \parbox{11cm}{for each $I\subset\mathbb N$ of $\card I=s$ and
       each subspace $F'\subseteq G$ of $\dim F'=l'$ there exists
       $x\in S_{F'}$ such that $\max_{n\in I}q_n(x)\le\delta$.}
  \end{equation}
  Choose an arbitrary vector $x_1\in S_G$ and a subset $I_1\subset\mathbb N$ of
  $\card I_1=s$ so that $\min_{n\in
    I_1}q_n(x_1)\ge\max_{n\in\mathbb N\setminus I_1}q_n(x_1)$. 
  It follows from~(\ref{qn_ge_r}) that there exists an index $n_1$
  such that $q_{n_1}(x_1)\ge r$; we can assume that $n_1\in I_1$.
  On the other hand, the definition of $s$ implies that
  $q_n(x_1)\le\delta$ whenever $n\notin I_1$. It
  follows from the definition of $L_m$ that there exists a subspace
  $G_{m-1}$ of $G$ of dimension $L_{m-1}$ so that
  \begin{equation}
    \label{max-max}
    \max\limits_{x\in S_{G_{m-1}}}\max\limits_{n\in I_1}q_n(x)\le
    \min\limits_{x\in S_{G_{m-1}}}\max\limits_{n\in I_1}q_n(x)+\delta\le
    2\delta,
  \end{equation}
  where the last inequality follows from~(\ref{s-set}). 
  
  Next, pick an $x_2\in S_{G_{m-1}}$ and $I_2\subset\mathbb N\setminus
  I_1$ so that $\card I_2=s$ and $\min_{n\in
    I_2}q_n(x_2)\ge\max_{n\notin I_1\cup I_2}q_n(x_2)$. Again, it
  follows from~(\ref{qn_ge_r}) that there exists an index $n_2$ such
  that $q_{n_2}(x_2)\ge r$; we can assume that $n_2\in I_1\cup I_2$.
  By~(\ref{max-max}), $q_n(x_2)\le 2\delta<r$ for each $n\in I_1$, so
  that $n_2\in I_2$.  Again, $q_n(x_2)\le\delta$ whenever $n\notin
  I_1\cup I_2$.  We can choose a subspace $G_{m-2}$ of $G_{m-1}$ of
  dimension $L_{m-2}$ so that
  $$\max\limits_{x\in S_{G_{m-2}}}\max\limits_{n\in
    I_2}q_n(x)\le2\delta.$$

  Proceeding this way, we obtain a sequence of vectors $x_1,\dots,x_m$
  and disjoint sets $I_1,\dots,I_m$ of cardinality $s$, and indices
  $n_1,\dots,n_m$ such that for each $i=1,\dots,m$ we have
  $n_i\in I_i$ and $q_{n_i}(x_i)\ge r$. Also,
  \begin{displaymath}
    q_n(x_i)\le
    \begin{cases}
      2\delta & \mbox{ if }n\in I_1\cup\dots\cup I_{i-1},\mbox{ and}\\
      \delta & \mbox{ if }n\notin I_1\cup\dots\cup I_i,
    \end{cases}
  \end{displaymath}
  hence $q_n(x_i)\le 2\delta$ whenever $n\notin I_i$. If $n\in I_i$
  then $q_n(x_i)\le\norm{x_i}=1$.

  Put $x=\sum_{i=1}^mx_i$, then for every $n\in\mathbb N$ we have
  $q_n(x)\le 1+m\cdot 2\delta\le 2$. On the other hand, 
  $$r\le q_{n_i}(x_i)\le q_{n_i}(x)+q_{n_i}(x-x_i)\le q_{n_i}(x)+2m\delta,$$
  so that $q_{n_i}(x)\ge
  r-2m\delta=\frac{r}{2}$
  for each $i=1,\dots,m$. It follows from the definition of $m$ that
  there can be at most $m-1$ indices $n$ such that
  $q_n(x)\ge\frac{r^2}{4}\norm{x}$, hence
  $\frac{r^2}{4}\norm{x}>\frac{r}{2}$. It follows that
  $\norm{x}>\frac{2}{r}$, so that $q_n(x)\le 2<r\norm{x}$ for every
  $n\in\mathbb N$, contradiction.
\end{proof}

Now we are ready to prove that $U$ is FSS.

\begin{theorem}
  The operator $U$ constructed in~(\ref{U-def}) is FSS for all $p\le
  2\le q$ unless $p=q=2$ or $p=1$ and $q=\infty$.
\end{theorem}

\begin{proof}
  In view of Theorem~\ref{U}(\ref{U-fss}) we may assume that $q=p'$.
  Recall that $U=\bigoplus_{n=1}^\infty U_n$ is composed of blocks
  $U_n\colon X_n\to Y_n$, where $X_n=\ell_p^{2^n}$ and
  $Y_n=\ell_q^{2^n}$.  For each $n$, let $P_n\colon\ell_p\to X_n$ be
  the canonical projection. For $x\in\ell_p$ put
  $q_n(x)=\bignorm{U_nP_nx}$. By Theorem~\ref{U}(\ref{norm-one}) we
  have $q_n(x)\le\norm{x}$.

  Assume that $U$ is not FSS. Then there exists a constant $C$ such
  that there are subspaces $G$ of $\ell_p$ of arbitrarily large dimension 
  such that the restriction of $U$ to $G$ is a $C$-isomorphism.
  Let $x\in S_G$, write
  $x=\sum_{n=1}^\infty x_n$ where $x_n\in X_n$,
  then $\norm{Ux}\ge\frac{1}{C}$. On the other hand,
  \begin{displaymath}
    \norm{Ux}^q=\sum\limits_{n=1}^\infty\norm{U_nx_n}^q\le
    \max\limits_{n\in\mathbb N}\norm{U_nx_n}^{q-p}
    \sum\limits_{n=1}^\infty\norm{U_nx_n}^p\le
    \max\limits_{n\in\mathbb N}q_n(x)^{q-p}.
  \end{displaymath}
  Hence, $\max\limits_{n\in\mathbb N}q_n(x)\ge C^{\frac{q}{p-q}}$.
  
  It follows from Lemma~\ref{flat-FDD} that there exists
  $\varepsilon>0$ such that for every $k$ and for every
  $G\subseteq\ell_p$ of sufficiently large dimension there exists a
  subspace $F$ of $G$ and an index $n$ such that $\dim F=k$ and
  $q_n(x)\ge\varepsilon$ for all $x\in S_F$. This implies that the
  restriction of $U_nP_n$ to $F$ is a
  $\frac{1}{\varepsilon}$-isomorphism. Put $E=P_n(F)$, then $E$ is a
  $k$\hyph dimensional subspace of $X_n$, and $U_n$ is a
  $\frac{1}{\varepsilon}$-isomorphism on~$E$.  In view of
  Theorem~\ref{K-convex} we may assume that $E$ is $2$-isomorphic to
  $\ell_2^k$ and $U_n(E)$ is $K$-complemented in $\ell_q^{2^n}$.
  
  Let $V_n$ be the canonical isometry between $L_p^N$ and
  $X_n=\ell_p^N$, where $N=2^n$, i.e., the isometry that maps
  $\chi_{[\frac{i-1}{N},\frac{i}{N}]}$ into $N^{-\frac{1}{p}}e_i$. It
  follows that
  $\norm{V_nx}_{\ell_r^N}=N^{\frac{1}{r}-\frac{1}{p}}\norm{x}_{L_r^N}$
  for every $x\in L_p^N$ and every $r\in[p,q]$. It follows from the
  definition of $U_n$ and Remark~\ref{H-norms} that
  \begin{displaymath}
    \bignorm{U_nV_n}_{L_2^N,\ell_2^N}=
      N^{\frac{1}{2}-\frac{1}{p}}\norm{U_n}_{\ell_2^N,\ell_2^N}
     =N^{\frac{1}{2}-\frac{1}{p}-\frac{1}{q}}\norm{H_n}_{\ell_2^N,\ell_2^N}=1
  \end{displaymath}
  and
  \begin{displaymath}
    \bignorm{U_nV_n}_{L_1^N,\ell_\infty^N}=
      N^{1-\frac{1}{p}}\norm{U_n}_{\ell_1^N,\ell_\infty^N}
     =N^{1-\frac{1}{p}-\frac{1}{q}}\norm{H_n}_{\ell_1^N,\ell_\infty^N}
     =1.
  \end{displaymath}
  Now applying Theorem~\ref{dimension} to $U_nV_n$ and $V_n^{-1}(E)$ we
  obtain a contradiction with the fact that $k=\dim E$ was chosen
  arbitrarily.
\end{proof}

\begin{remark}
  If $p=q=2$ then $U$ is an isometry, hence not FSS. Consider the
  case when $p=1$ and $q=\infty$.  The preceding proof doesn't work
  since now we cannot use Theorem~\ref{K-convex}.  Actually, $U$ is
  not FSS in this case. Indeed, we now have $U_n=H_n$, and it is easy
  to see that one finds all the $n$ Rademacher sequences of length
  $2^n$ among the columns of $H_n$.  That is, there are $1\le
  j_1,\dots,j_n\le n$ such that $H_ne_{j_i}=r_i$, where $r_i$ is the
  $i$-the Rademacher sequence (of length $2^n$). Since $r_1,\dots,r_n$
  span an isometric copy of $\ell_1^n$ in $\ell_\infty^N$, it follows
  that the restriction of $H_n$ to the span of $e_{j_1},\dots,e_{j_n}$
  preserves a copy of $\ell_1^n$.
\end{remark}

\begin{question}
  Are there any other closed ideals in $\Lpq$? In view of the diagram
  at the beginning of our paper this question can be subdivided in the
  following subquestions:
  \begin{enumerate}
    \item Is $\JIpq$ equal to $\fss\cap\J2$?
       If not, is $\fss\cap\J2$ an immediate successor of 
        $\JIpq$?
    \item Is $\fss$ an immediate successor of $\fss\cap\J2$? 
        More generally, are there any immediate successors of    
        $\fss\cap\J2$, other than $\J2$?
    \item Is $\fss\vee\J2$ an immediate successor of $\J2$?
    \item Is $\fss\vee\J2$  equal to $\Lpq$?  
 \end{enumerate}
\end{question}

\begin{question}
  Suppose again that $U$ is the operator defined in~(\ref{U-def}).
  Since $U$ is FSS, we have $\iJ^U\subseteq\fss$.  Does $\iJ^U$ equal
  $\fss$?
\end{question}


\end{document}